\documentclass[12pt,a4paper]{amsart}
\usepackage{graphicx}

\usepackage{tikz}
\usetikzlibrary{positioning}

\newtheorem{thm}{Theorem}

\title[Minimal Grid Diagrams of the Prime Alternating Knots]%
{Minimal Grid Diagrams of the \\ Prime Alternating Knots with 13 Crossings}

\author[H. J. Lee]{Hwa Jeong Lee}
\address{Hwa Jeong Lee\\ Dongguk University WISE\\
123, Dongdae-ro, Gyeongju-si, Gyeongsangbuk-do 38066, Republic of Korea}
\email{hjwith@dongguk.ac.kr}

\author[A. Stoimenow]{Alexander Stoimenow}
\address{Alexander Stoimenow\\ Dongguk University WISE\\
123, Dongdae-ro, Gyeongju-si, Gyeongsangbuk-do 38066, Republic of Korea}
\email{stoimeno@stoimenov.net}

\author[G. T. Jin]{Gyo Taek Jin}
\address{Gyo Taek Jin\\ Korea Advanced Institute of Science and Technology\\
291 Daehak-ro, Yuseong-gu, Daejeon 34141, Republic of Korea}
\email{trefoil@kaist.ac.kr}

\def\h#1#2#3{\draw[black,line width=0.75pt, cap=round] (#1,#2)--(#1+#3,#2);}
\def\w{0.2}
\def\v#1#2#3{
       \draw[white,line width=2.5pt] (#1,#2+\w)--(#1,#2+#3-\w);
       \draw[black,line width=0.75pt, cap=round] (#1,#2)--(#1,#2+#3);}
\def\knum#1#2#3#4{
\draw (-0.5,-0.5) (15.0,17);
\filldraw[black] (0,16.5) node[anchor=west]{$#1#2#3$}; 
}

\begin{document}

\begin{abstract}
A \emph{knot} is a closed loop in space without self-inter\-section. Two knots are equivalent if there is a self homeomorphism of space bringing one onto the other. An \emph{arc presentation} is an embedding of a knot in the union of finitely many half planes with a common boundary line such that each half plane contains a simple arc of the knot. The minimal number of such half planes among all arc presentations of a given knot is called the \emph{arc index} of the knot. A knot is usually presented as a planar diagram with finitely many crossings of two strands where one of the  strands goes over the other. 
A \emph{grid diagram} is a planar diagram which is a non-simple rectilinear
polygon such that vertical edges always cross over horizontal edges at all
crossings.
It is easily seen that an  arc presentation gives rise to a grid diagram and vice versa. It is known that the arc index of an alternating knot is two plus its minimal crossing number. There are 4878 prime alternating knots with minimal crossing number 13. We obtained minimal arc presentations of them in the form of grid diagrams having 15 vertical segments. 
This is a continuation of the works on prime alternating knots of 11 crossings and 12 crossings.
\end{abstract}
\date{\today}
\makeatletter
\@namedef{subjclassname@2020}{\textup{2020} Mathematics Subject Classification}
\makeatother

\maketitle

\section{Introduction}
A \emph{grid diagram\/} is a diagram of a knot or a link with finitely many vertical segments and the same number of horizontal segments in which the vertical segments cross over the horizontal segments at all crossings.
Given a grid diagram of a knot, consider pulling middle points of the horizontal segments back to the diagram so that they touch a vertical axis. Then each half plane determined by the axis and a vertical segment contains an arc of the knot. Such a position of a knot is called an \emph{arc presentation}. 
The \emph{arc index} of a knot $K$, denoted $\alpha(K)$, is the minimum number of half planes among all arc presentations of $K$, or equivalently, the minimum number of vertical segments among all grid diagrams of $K$~\cite{C1995}.

\begin{figure}[h!]
\begin{tikzpicture}[scale=0.25]
\input 13a0001t.tex
\end{tikzpicture}
\caption{A minimal grid diagram of the knot $13a1$}\label{fig:13a1grid}
\end{figure}
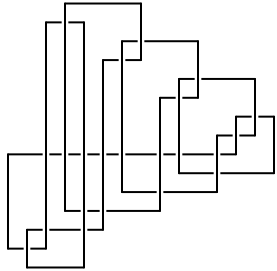

There are 4878 prime alternating knots with 13 crossings~\cite{knotinfo}.
We wrote a series of  programs on Maple to find their minimal grid diagrams.
According to Theorem\,\ref{thm:Bae-Park}, the minimal grid diagrams of the 13 crossing knots have  15 vertical segments.
\begin{thm}[Bae-Park \cite{BP2000}]\label{thm:Bae-Park}
If $K$ is an alternating knot, then the arc index of $K$ is the minimal crossing number plus two.
\end{thm}

 For each of the mentioned 4878  knots, we pursued the following steps:
\begin{enumerate}
\item Using the DT notation from  Knotscape~\cite{knotscape}, the regions divided by the knot diagram are described.
\item Considering the knot diagram as a planar graph, a spanning tree is generated whose contraction leads to an arc presentation of the knot.
\item The arc presentation  is converted to a grid diagram in various forms such as 2d and 3d graphics, latex picture commands and a sequence of 3d coordinates of vertices as a polygonal knot.
\item Using Knotplot~\cite{knotplot}, a DT notation for the grid diagram is obtained.
\item Using Knotscape~\cite{knotscape}, the grid diagram is confirmed to be the same as the original knot up to taking mirror images.
\end{enumerate}

\section{An example}
\def\mesh{
\color{cyan}
\multiput(0,0)(0,20){10}{\line(1,0){180}}
\multiput(0,0)(20,0){10}{\line(0,1){180}}
\color{red}
\put(35,5){40}\put(75,5){80}\put(115,5){120}\put(155,5){160}
\put(-10,38){40}\put(-10,78){80}\put(-15,118){120}\put(-15,158){160}}
\begin{figure}[h!]\tiny
\begin{picture}(180,180)
\put(0,0){\includegraphics[width=0.50\textwidth]{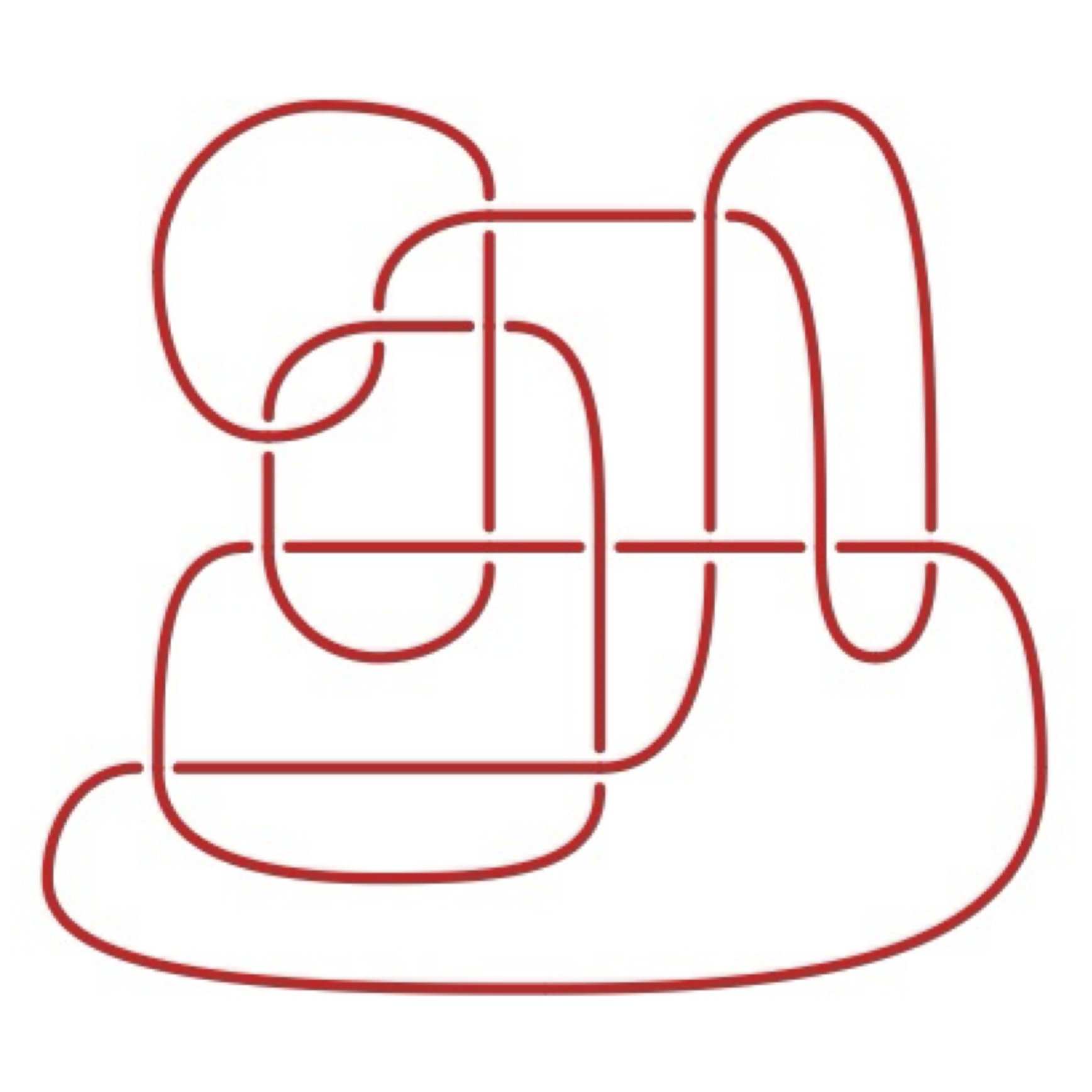}}%
\color{black}
\put(138,118){\boldmath$\Downarrow$}
\put(128,118){1} \put(142,75){2} \put(155,122){3} \put(111,115){4}
\put(119,65){5}  \put(50,57){6}  \put(91,20){7} \put(142,92){8} 
\put(125,92){9}  \put(105,92){10} \put(84,92){11} \put(59,82){12}
\put(28,65){13} \put(65,38){14} \put(88,65){15} \put(93,124){16} 
\put(66,117){17} \put(43,127){18} \put(32,97){19}  \put(59,63){20} 
\put(70,100){21} \put(83,132){22}  \put(47,152){23}  \put(52,102){24} 
 \put(58,140){25}  \put(93,148){26} 
\end{picture}
\caption{$13a1$}\label{fig:13a1}
\end{figure}

\begin{figure}[t!]\tiny
\centerline{
\begin{picture}(180,180)
\put(0,0){\includegraphics[width=180pt]{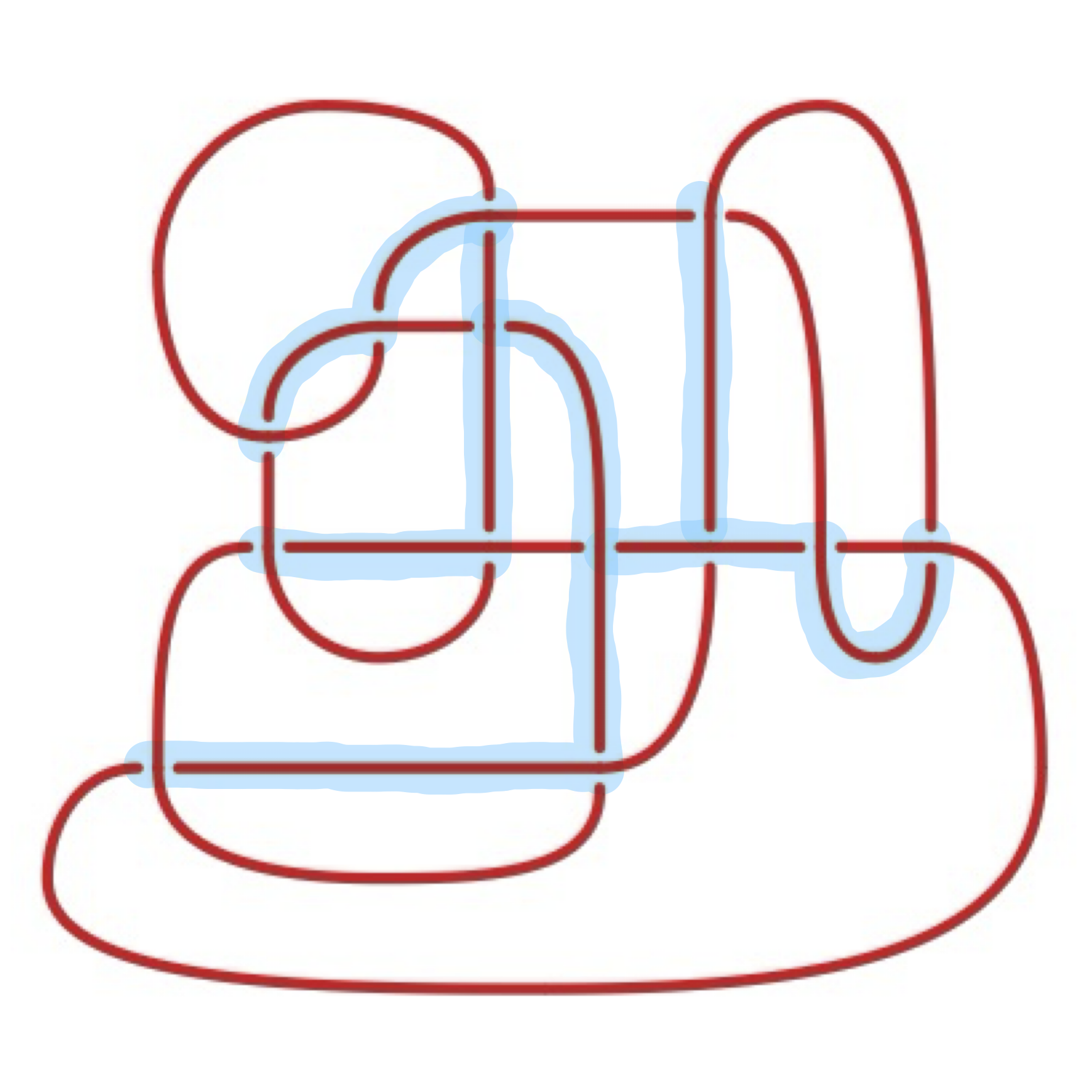}}
\put(90,-5){\small(a)}
\color{blue}\bf
\put(81,144){\circle{10}}
\put(66,143){1}\put(82,133){2}\put(47,125){3}\put(83,108){4}
\put(96,122){5}\put(60,92){6}\put(92,68){7}\put(106,81){8}
\put(60,56){9}\put(121,81){10}\put(106,115){11}\put(140,63){12}
\put(152,140){\boldmath$\star$}
\end{picture}
\begin{picture}(180,180)
\put(0,0){\includegraphics[width=180pt]{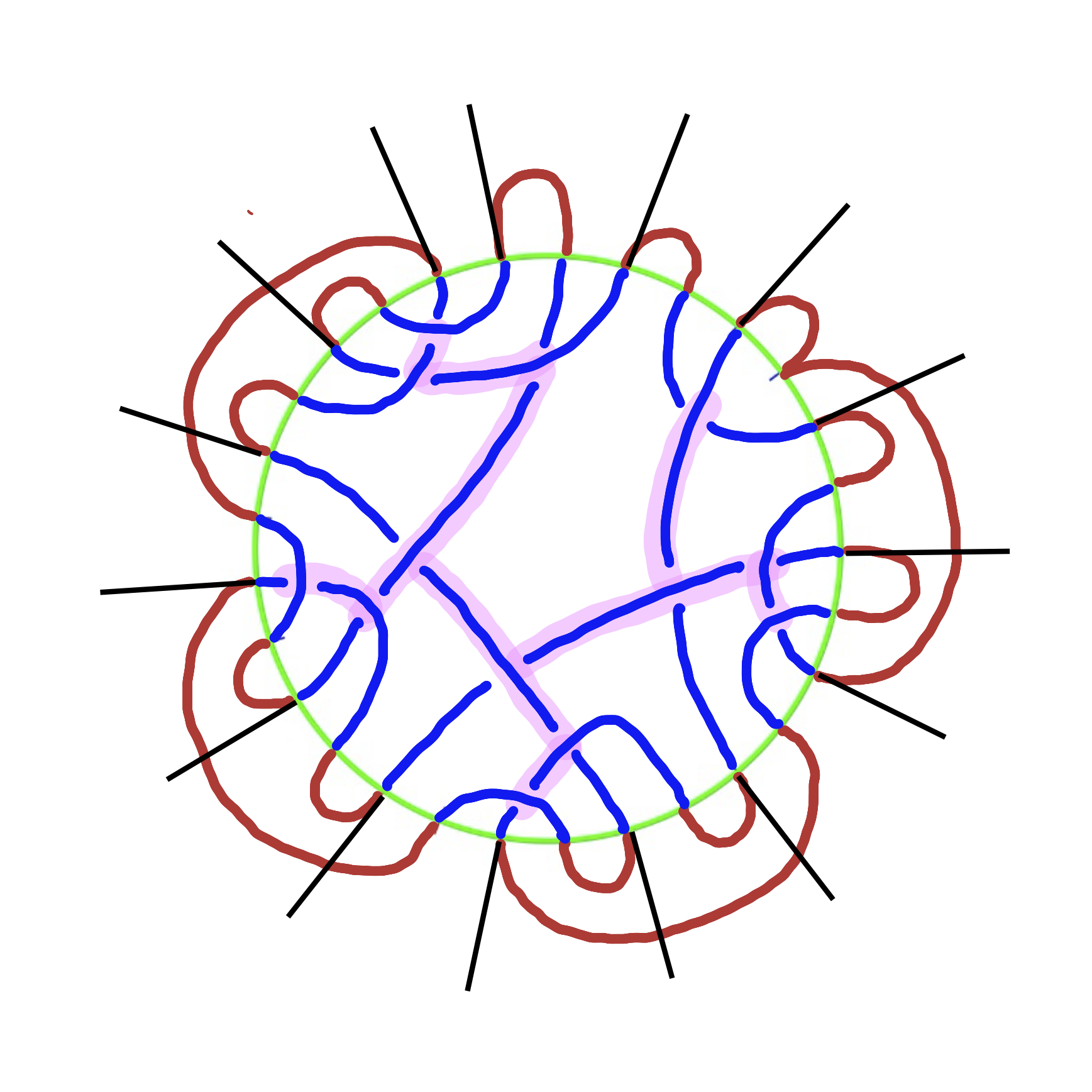}}
\put(90,-5){\small(b)}\bf
\put(88,120){\circle{7}}
\put(124,115){1}\put(104,115){2}\put(104,93){3}\put(94,80){4}
\put(80,80){5}\put(80,95){6}\put(94,115){7}\put(60,105){8}
\put(83,127){9}\put(64,67){10}\put(48,90){11}\put(96,64){12}
\put(73,51){13}\put(119,95){14}\put(114,74){15}
\put(103,163){$(2,7)$}\put(141,144){$(1,3)$}\put(161,120){$(2,14)$}
\put(161,92){$(4,15)$}\put(154,50){$(1,14)$}\put(138,25){$(3,12)$}
\put(104,10){$(5,13)$}\put(64,10){$(12,15)$}\put(38,20){$(4,10)$}
\put(11,42){$(6,11)$}\put(5,87){$(10,13)$}\put(7,116){$(5,8)$}
\put(25,142){$(7,9)$}\put(46,161){$(8,11)$}\put(75,165){$(6,9)$}
\put(135,121){$\star$}
\put(95,155){\line(1,5){5}}\put(100,180){\vector(4,-1){20}}
\end{picture}
}
\caption{A  spanning tree on a diagram of $13a1$ and its corresponding spokes}\label{fig:13a1_Tree}
\end{figure}
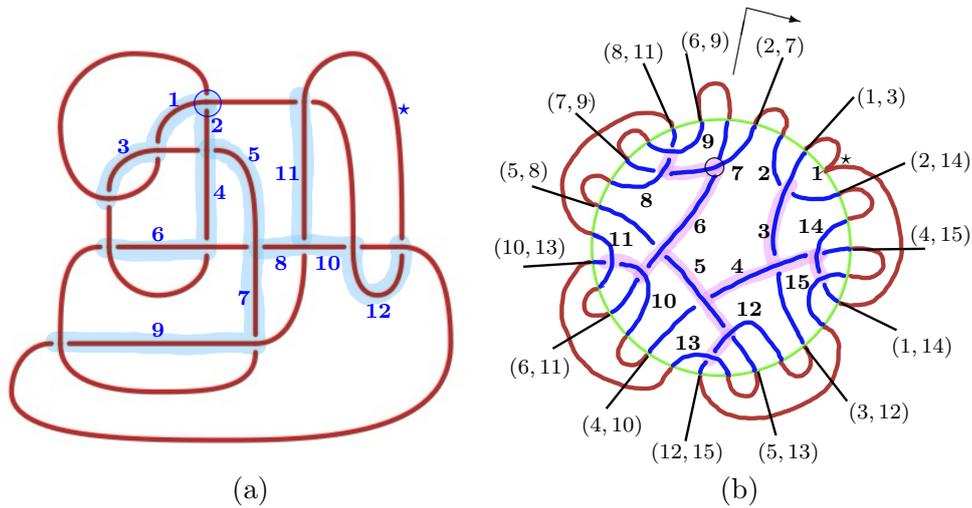

We give a  detailed description of the steps (1)--(3) using the knot $13a1$.
The diagram in Figure\,\ref{fig:13a1}  shows a minimal diagram of the knot $13a1$ obtained from \cite{knotinfo} with the edge labels
compatible with the DT notation 	[4, 8, 10, 14, 2, 16, 20, 6, 22, 24, 12, 26, 18]\footnote{from Knotscape~\cite{knotscape}} in the following sense. 
At a crossing incident to the four edges with lables $i$, $i+1$, $j$, $j+1$, modulo $26$,  we may assume that $i$ and $j+1$ are odd and $i+1$ and $j$ are even. 
Then $i+1$ appears in the $(j/2+1)$-st (modulo 13) place in the DT notation. We orient the knot in the direction as the labels increase.
The diagram divides the plane into 15 regions, namely,
$$\begin{aligned}
&[2, 8],[6, 14],[12, 20],[18,24],[-1, 4, -9],[1, -8, 3], [-5,10,-15],\\
&[-11,16,-21],[-17,22,-25],[-18,25,23],[-4,-26,-22,-16,-10], \\
&[-6,15,11,-20,13],[-12,21,17,-24,19],[-2,9,5,-14,7],\\
&[-3,-7,-13,-19,-23,26]
\end{aligned}$$
where each sequence describes the oriented edges on the boundary. The numbers with minus sign indicate that the corresponding edges are oriented in the opposite of the orientation of the knot. All  regions except the only unbounded region $[-3,-7,-13,-19,-23,26]$ are oriented counterclockwise. 

The thickened edges of Figure\,\ref{fig:13a1_Tree}(a) form a spanning tree of the diagram. This tree is rooted at  the crossing between the edges 25 and 26.  From the root,  it grows as the following sequence $e_i, i=1,\ldots,12$, of oriented edges  indicates.
\begin{equation*}\label{eq:tree}
-25,  -22,  18,  -21,  -16,  12, -15, -10,  6,  -9, -4,  2.        
\end{equation*}
For each $i=1,\ldots,12$,  edge $e_i$ of the tree is labeled {\boldmath$ i$} in Figure\,\ref{fig:13a1_Tree}(a). Let $c$ denote the crossing number. The edges of the spanning tree are chosen to satisfy the following conditions:
\begin{itemize}
\item For each $j=1,\ldots,c-1$, $e_1\cup\cdots \cup e_j$ is a tree.
\item For each $j=1,\ldots,c-1$, $e_1\cup\cdots \cup e_j$ does not separate untouched crossings into two parts.
\item For each $j=1,\ldots,c-2$, $e_1\cup\cdots \cup e_j\cup e$ has no loop where $e$ is the one-step extension of the oriented edge $e_j$.
\end{itemize}

The proof of Theorem~\ref{thm:Bae-Park} shows that for any prime link diagram a spanning tree can be chosen to satisfy the above three conditions~\cite{BP2000,Jin2012}.

The knot diagram in Figure\,\ref{fig:13a1_Tree}(b) is obtained from that of Figure\,\ref{fig:13a1_Tree}(a) by a plane isotopy, where the circle is the boundary of a tubular neighborhood of the spanning tree.
There are 14 edges outside a tubular neighborhood in (a). The edge in (a) labeled $\star$ is broken in the middle so that there are 15 arcs outside the circle in (b).

Inside the circle in (b), there are 14 arcs numbered 2 to 15. We may put each of these arcs in a horizontal plane at the height given by its number.  These heights are chosen in the following order:
\begin{itemize}
\item At the root, we give 7 to $e_1$ and  6 to  $e_2$. 
\item At the end of $e_1$, we give 8 to $e_3$. 
\item At the end of  $e_2$, we give 5 to $e_5$. 
\item At the end of  $e_3$, we give 9 to the overcrossing arc. 
\item $e_4$ inherits 6 from $e_2$. 
\item At the end of  $e_4$, we give 10 to $e_6$. 
\item At the end of  $e_5$, we give 4 to $e_8$. 
\item At the end of  $e_6$, we give 11 to  the overcrossing arc. 
\item $e_7$ inherits 5 from $e_5$. 
\item At the end of  $e_7$, we give 12 to  $e_9$. 
\item At the end of  $e_8$, we give 3 to  $e_{11}$. 
\item At the end of  $e_9$, we give 13  to the overcrossing arc. 
\item $e_{10}$ inherits 4 from $e_8$. 
\item At the end of  $e_{10}$, we give 14 to  $e_{12}$. 
\item At the end of  $e_{11}$, we give 2 to  the undercrossing arc. 
\item At the end of  $e_{12}$, we give 15 to  the overcrossing arc. 
\item Finally, the middle point of the the extension of $e_{12}$ is given the height 1.
\end{itemize}

Now we need to put the 15 arcs outside the circle into vertical half planes. The vertical positions of their endpoints are already determined by the heights of the arcs inside the circle. As the part inside the circle is shrunk to points on an axis, the eleven innermost arcs can be easily tilted to lie in  vertical half planes on one of their ends.
The remaining four arcs have one of their ends interleaved in the ends of their inner arcs. In order not to make crossing changes, their non-interleaved ends must be tilted toward the interleaved ones.
The spokes in Figure\,\ref{fig:13a1_Tree}(b) were thus created. We read the vertical intervals of their endpoints as the broken arrow directs:
\begin{equation*}\label{eq:12a1wheel}
\begin{aligned}
&[2, 7], [1, 3], [2, 14], [4, 15], [1, 14], [3, 12], [5, 13], [12, 15], [4, 10], \\
&[6, 11], [10, 13], [5, 8], [7, 9], [8, 11], [6, 9]
\end{aligned} 
\end{equation*}

This is the sequence of intervals of the vertical segments in the grid diagram in Figure\,\ref{fig:13a1grid}.

The idea behind this construction is originally due to Bae and Park~\cite{BP2000}, modified by two authors using the spanning tree~\cite{Jin2012}.
Interested readers are encouraged to read these articles.

\section{Shaping the grid diagrams}
In order for the 4878 grid diagrams to look less random, we modify them using the following moves:
\begin{enumerate}
\item A set of consecutive parallel segments on one end is moved to the other end.
\item Rotation about the diagonal or antidiagonal axis.
\item Reflection about a horizontal or vertical line.
\end{enumerate}

\begin{figure}[h]
\leftline{
\begin{tikzpicture}[scale=0.23]
\input 13a0001t.tex
       \draw[green,line width=0.75pt] (1,13.5)--(10,13.5)
                                                     --(10,15.5)--(1,15.5)--(1,13.5);
            \draw[white,line width=0.75pt] (0.5,0.5) (15.5,13.5);                                               
\end{tikzpicture}
\begin{tikzpicture}[scale=0.23]
\filldraw[black] (2,7.5) node[anchor=west]{$\overset{(1)\;}\Longrightarrow$}; 
            \draw[white,line width=0.75pt] (0.5,0.5) (6.5,13.5);                                               
\end{tikzpicture}
\begin{tikzpicture}[scale=0.23]
\input 13a0001t2.tex
       \draw[green,line width=0.75pt] (0.5,-1.5)--(5.5,-1.5)
                                                     --(5.5,7.5)--(0.5,7.5)--(0.5,-1.5);
      \draw[white,line width=0.75pt] (0.5,-0.5) (15.5,13.5);
\end{tikzpicture}}

\medskip
\rightline{
\begin{tikzpicture}[scale=0.23]
\filldraw[black] (-1,7.5) node[anchor=west]{$\overset{(1)\;}\Longrightarrow$}; 
            \draw[white,line width=0.75pt] (0.5,0.5) (3.5,13.5);                                               
\end{tikzpicture}
\begin{tikzpicture}[scale=0.23]
\input 13a0001t3.tex
       \draw[white,line width=0.75pt] (5.5,-1.5) (20.5,13.5);
       \draw[green,line width=0.75pt] (6.5,-0.5)--(19.5,12.5);       
\end{tikzpicture}
\begin{tikzpicture}[scale=0.23]
\filldraw[black] (2,7.5) node[anchor=west]{$\overset{(2)\;}\Longrightarrow$}; 
            \draw[white,line width=0.75pt] (0.5,0.5) (6.5,13.5);                                               
\end{tikzpicture}
\begin{tikzpicture}[scale=0.23]
\input 13a0001u.tex
      \draw[white,line width=0.75pt] (0.5,0.5) (15.5,15.5);
\end{tikzpicture}}
\caption{Shaping a grid diagram of $13a1$}
\end{figure}
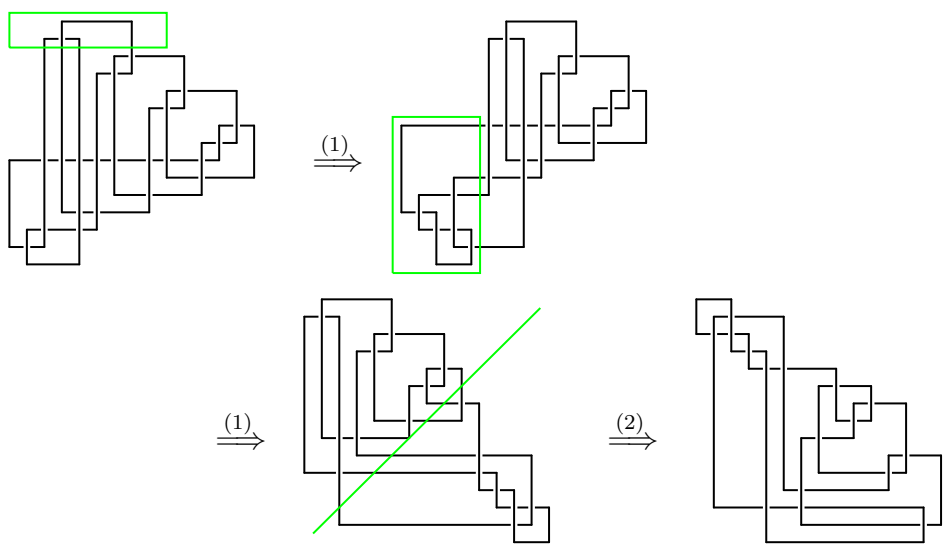

\section{Minimal grid diagrams of 13 crossing prime alternating knots}

\def\h#1#2#3{\draw[black,line width=0.6pt, cap=round] (#1,#2)--(#1+#3,#2);}
\def\w{0.2}
\def\v#1#2#3{
       \draw[white,line width=2.0pt] (#1,#2+\w)--(#1,#2+#3-\w);
       \draw[black,line width=0.6pt, cap=round] (#1,#2)--(#1,#2+#3);}
\def\knum#1#2#3#4{
\draw (-0.5,-0.5) (15.0,17);
\filldraw[black] (0,16.5) node[anchor=west]{$#1#2#3$}; 
}
  
{\scriptsize
\noindent
\input 13a4878ejc.tex
\begin{tikzpicture}[scale=0.115]  \draw (-0.5,-0.5) (15.0,0.5);
\end{tikzpicture}  
\begin{tikzpicture}[scale=0.115]  \draw (-0.5,-0.5) (15.0,0.5);
\end{tikzpicture}  
\begin{tikzpicture}[scale=0.115]  \draw (-0.5,-0.5) (15.0,0.5);
\end{tikzpicture}  
\begin{tikzpicture}[scale=0.115]  \draw (-0.5,-0.5) (15.0,0.5);
\end{tikzpicture}  
\begin{tikzpicture}[scale=0.115]  \draw (-0.5,-0.5) (15.0,0.5);
\end{tikzpicture}  
}
\section*{Acknowledgments}
  This work was supported by the National Research Foundation of Korea(NRF) grant funded by the Korea government(MSIT) (No. 2023 R1A2C1003749).

\end{document}